\documentstyle[amsthm,12pt,pb-diagram]{article}

\newcommand{\ul}[1]{\underline{#1}}
\newcommand{\RA}{\Rightarrow}

\newcommand{\bss}{bar spectral sequence}
\newcommand{\etwo}{E^2_{*,*}}

\newcommand{\ko}[1]{\underline{KO}_{\, #1}}
\newcommand{\ku}[1]{\underline{KU}_{\, #1}}
\newcommand{\hgy}[1]{H_*{#1}}
\newcommand{\tor}[1]{\mbox{Tor}^{{#1}}_{*,*}(\ints/2, \ints/2)}
\newcommand{\torh}[1]{\mbox{Tor}^{H_*{#1}}_{*,*}(\ints/2, \ints/2)}
\newcommand{\hko}[1]{H_*\underline{KO}_{\, #1}}
\newcommand{\hku}[1]{H_*\underline{KU}_{\, #1}}
\newcommand{\olz}[1]{\overline{z}_{#1}}
\newcommand{\es}[1]{e^{\circ #1}}

\newcommand{\be}{\beta}
\newcommand{\la}{\lambda}
\newcommand{\lai}{\lambda^{-1}}
\newcommand{\nui}{\nu^{-1}}
\newcommand{\Ga}{\Gamma}
\newcommand{\ga}{\gamma_{2^j}}
\newcommand{\si}{\sigma}

\newcommand{\ints}{\mbox{\bf Z}}
\newcommand{\rls}{\mbox{\bf R}}
\newcommand{\cmpx}{\mbox{\bf C}}

\newcommand{\spec}[1]{\ul{E}_{\, #1}}

\begin{document}

\title{
        The Hopf Rings for $KO$ and $KU$
}
\author{
        Dena S. Cowen Morton$^{\mbox{a}}$,  Neil Strickland $^{\mbox{b}}$
}

\maketitle

\begin{center}
$^{\mbox{a}}$ {\it Department of Mathematics, Xavier University, Cincinnati, OH 45207,  USA}

$^{\mbox{b}}$ {\it Department of Pure Mathematics, University of Sheffield, 
  Sheffield S3 7RH, UK}

\end{center}

\begin{abstract}
 In this work we compute the Hopf rings for the spectra representing
 orthogonal and unitary $K$-theory, $KO$ and $KU$.  Specifically, we
 use Hopf ring properties to find and simplify the ordinary mod 2
 homology groups for these spectra.  Our main tool is the \bss, which
 allows us to advance from one space to the next in the spectrum.

\


\end{abstract}

\theoremstyle{plain}
\newtheorem{lem}{Lemma}[section]
\newtheorem{prop}[lem]{Proposition}

\section*{Contents}

\noindent 1. Introduction

\noindent 2. Hopf Rings 

\noindent 3. The Homotopy Elements $[x]$

\noindent 4. Properties for $\hko{0}=\hgy{(\ints \times BO)}$

\noindent 5. Machinery

\noindent 6. The Calculation of the Hopf Ring for $\hko{*}$

\noindent 7. Properties for $\hku{0}=\hgy{(\ints \times BU)}$

\noindent 8. The Calculation of the Hopf Ring for $\hku{*}$

\section{Introduction}

The object of this paper is to compute the mod 2 Hopf rings,
$H_*(\ko{*};\ints/2)$ and $H_*(\ku{*};\ints/2)$, for the periodic Bott
spectra $KO$ and $KU$. (For the remainder of this paper, we use the
notation $\hgy{X}$ to mean $H_*(X;\ints/2)$.)  The spaces in these
spectra are the infinite classical groups and their coset spaces, and
their homology was first calculated in ~\cite{cartan}, but the Hopf ring
structure was first determined in the second author's unpublished
thesis ~\cite{strick}.  The presentation given here serves as an
introduction to the first author's much more intricate work on the
connective spectrum $bo$.  The Hopf ring viewpoint turns out
to be very convenient for understanding the homological effect of
various maps between classical groups and fibrations of their
connective covers, for example in transferring the results
of ~\cite{hoanst} from a complex context to a real context; details will
be given elsewhere.

We next give a brief statement of our results (using some standard
Hopf ring notation that will be recalled in
section ~\ref{sec-hopf-rings}).  We push forward the usual generators
of $\hgy{\rls P^\infty}$ along the inclusion
$\rls P^\infty=BO(1)\subseteq\{1\}\times BO\subseteq \ints\times
BO=\ko{0}$ to get elements $z_k\in H_k\ko{0}$.  We also write $z_k$
for the image of this element in $H_k\ku{0}$ (which is zero when $k$
is odd).  We put $\olz{k}=z_k/z_0=z_k*z_0^{-1}=z_k*[-1]$.  It is well-known that
        $$KO_{\,*}= {\ints}[\eta,  \beta,\lambda^{\pm 1}]/
           (\eta^{3},2 \eta,  \eta\beta, \beta^{2} - 4 \lambda), $$
where $\deg(\eta) = 1,  \deg(\beta)=4,$ and $\deg(\lambda)=8,$ and that 
        $$KU_{\,*}=\ints[\nu^{\pm 1}],$$ 
where $\deg(\nu)=2$.  This gives the ring-ring  elements
        $$[\eta]\in H_{0}\ko{-1}, [\be]\in H_{0}\ko{-4}, 
        [\la]\in H_{0}\ko{-8},\ \mbox{and} \ [\nu]\in H_{0}\ku{-2}.$$
(For more information about ring-ring elements of the form above, see section 3.)

We also have the element $[1]=z_0\in H_0\ko{0}$.

For each of the spaces $X$ under consideration, the ring $\hgy{X}$ is
either polynomial or exterior on countably many generators, possibly
with a polynomial generator inverted.  The generators are indexed by a
parameter $i$, which always runs from $0$ to $\infty$.  If the
generator to be inverted has the form $[x]$, then the inverse is
$[-x]$.  The map $x\mapsto [\la]\circ x$ gives an isomorphism
$H_*\ko{n}\rightarrow H_*\ko{n-8}$, so we need only describe
$H_*\ko{n}$ for $0\leq n<8$.  Similarly, we have $H_*\ku{n}\simeq
H_*\ku{n-2}$.

Our detailed answer for $KO$ is as follows.

\begin{tabular} {r l l }
        $\hko{0}=$ & $\hgy{(\ints\times BO)}$ & $=P(\olz{i}, [-1])$\\
        $\hko{1}=$ & $\hgy{(U/O)}$ & $=P(e\circ z_{2i})$\\ 
        $\hko{2}=$ & $\hgy{(Sp/U)}$ & $=P(\es{2}\circ z_{4i})$\\ 
        $\hko{3}=$ & $\hgy{(Sp)}$ & $=E(\es{3}\circ z_{4i})$\\ 
        $\hko{4}=$ & $\hgy{(\ints\times BSp)}$ &
                     $=P(\olz{4i}\circ[\be\lai], [-\be\lai])$\\ 
        $\hko{5}=$ & $\hgy{(U/Sp)}$ & $=E(e\circ z_{4i}\circ[\be\lai])$\\ 
        $\hko{6}=$ & $\hgy{(O/U)}$ & $=E(\olz{2i}\circ [\eta^2\lai])$\\ 
        $\hko{7}=$ & $\hgy{(O)}$ & $=E(\olz{i}\circ[\eta\lai])$\\ 
        $\hko{8}=$ & $\hgy{(\ints\times BO)}$ &
                     $=P(\olz{i}\circ[\lai], [-\lai])\cong \hko{0}$\\
\end{tabular}

\noindent The mod 2 homology for the first term, $\ko{0}=\ints \times BO$, 
is well known; accordingly, we will take $\ko{0}$ to be our starting point, and 
deloop to $\ko{1}$ using the \bss.  This process will continue, from 
$\ko{1}$ to $\ko{2}$, etc., throughout the 8-space cycle until we end with 
$\ko{8}=\ko{0}$.  We will prove that our spectral sequences collapse
and solve the extension problems using Hopf ring relations in
conjunction with the Frobenius and Verschiebung maps.

We also record the answer for $KU$; this is essentially well-known,
but it is convenient to have it stated in a way that allows easy
comparison with the case of $KO$.
\begin{tabular} {r l l }
        $\hku{0}=$ & $\hgy{(\ints\times BU)}$ & $=P(z_{2i}, [-1])$\\
        $\hku{1}=$ & $\hgy{(U)}$ & $=E(e\circ z_{2i})$\\ 
        $\hku{2}=$ & $\hgy{(\ints\times BU)}$ & 
        $=P(z_{2i}\circ[\nui],[-\nui])$\\ 
\end{tabular}

The Hopf ring relations that we need are as follows:
\begin{enumerate}
 \item $z_j\circ z_k=\frac{(j+k)!}{j! k!}z_{j+k}$
 \item $e^2 = e\circ z_1$
 \item $(\es{2})^2 = \es{2}\circ z_2$
 \item $(\es{3})^2 = 0$
 \item $\es{4}\circ[\la]=[\be]\circ\olz{4}$
 \item $e\circ[\eta]=\olz{1}$
 \item $\es{2}\circ[\be]=\olz{2}\circ[\eta^2]$
 \item $\olz{1}\circ[\be]=\olz{2}\circ[\be]=0$
 \item $\olz{1}\circ\olz{2i+1}=\olz{1}^2\circ \olz{2i}$
\end{enumerate}
Proofs are distributed throughout the paper.

We conclude our introduction with a brief discussion of various maps of spectra. 
 There
is a complexification map $m:KO\rightarrow KU$, a complex conjugation
map $c: KU\rightarrow KU$, and a map $f: KU\rightarrow KO$ that
forgets the complex structure.  These satisfy
\begin{enumerate}
 \item $cm=m$, $fc=f$, $c^2=1$, $fm=2$, $mf=1+c$
 \item $m_*(\eta)=0$, $m_*(\be)=2\nu^2$, $m_*(\la)=\nu^4$
 \item $c_*(\nu)=-\nu$
 \item $f_*(1)=2$, $f_*(\nu)=\eta^2$, $f_*(\nu^2)=\be$,
  $f_*(\nu^3)=0$.
\end{enumerate}
The maps $m$ and $c$ are ring maps and thus induce maps
of Hopf rings.  The map $f$ is a $KO$-module map, so it satisfies
$f_*(b*c)=f_*(b)*f_*(c)$ and $f_*(m_*(a)\circ b)=a\circ f_*(b)$ (for
$a\in H_*\ko{*}$ and $b,c\in H_*\ku{*}$).  Using these properties, one
can determine the effects of $m_*$, $f_*$ and $c_*$ on all the elements
of our Hopf rings.

\section{Hopf rings}
\label{sec-hopf-rings}

Let $R$ be a graded associative commutative ring with unit and 
let $CoAlg_{R}$ denote the category of graded cocommutative coassociative 
coalgebras with counit over $R$. Then a Hopf ring is a graded ring object 
in the category $CoAlg_{R}$. A Hopf ring includes a coproduct $\psi$, 
two products - the $*$-product and the $\circ$-product, conjugation $\chi$, and
relationships interlocking each of these maps. 

The primary example of Hopf rings is $F_*\ul{E}_*$, where $F_*$ is a
multiplicative homology theory and $\ul{E}_*$ is an $\Omega$-spectrum.
For more information about maps and properties of Hopf rings, please
see \cite{ravwil}. 

\section{The Homotopy Elements $[x]$}

The computations of $\hko{0}$ and $\hku{0}$ will require the use of
elements from the homotopy groups of $KO$ and $KU$, as ring-ring
elements.  In particular, we will make heavy use of the homotopy elements
$$[\eta],[\be],[\la], \mbox{and} \ [\nu].$$   More information is available in 
\cite{wil} and \cite{strick}.

Following \cite{ravwil}, suppose ${\cal S}$ is a homotopy category of
topological spaces (with certain properties), and $F_{*}(-)$ is an
associative commutative multiplicative unreduced generalized homology
theory with unit defined on ${\cal S}$.  If we let $G^{*}(-)$ be a
similar cohomology theory (also defined on ${\cal S}$), then
$G^{*}(-)$ has a representing $\Omega$-spectrum
$$\underline{E}_{\,*} =\{\underline{E}_{\,n}\}_{n \in {\bf Z}} \in
{\cal \mbox{gr}S},$$ 
i.e. $G^{n}(X) \simeq [X, \underline{E}_{\,n}]$
and $\Omega \ul{E}_{\,n+1} \simeq \ul{E}_{\,n}$ (with ${\cal
  \mbox{gr}S}$ the category of graded objects of $S$). Denoting the
two coefficient rings by $F_{*}$ and $G^{*}$, we let $x \in G^{n}$
have degree $-n$ in the coefficient ring.  Then $x \in G^{n} \simeq
[\mbox{point}, \underline{E}_{\,n}]$ and so we have a map in homology
$x_{*}:F_{*}\rightarrow F_{*}\ul{E}_{\,n} $.  We define $[x]\in
F_{0}\ul{E}_{\,n}$ to be the image of $1\in F_{*}$ under this map.

If we take $z \in G^{m}$ and $x, y \in G^{n}$,  then 
\begin{enumerate}
 \item $[z] \circ [x] = [zx]=[-1]^{\circ nk}\circ [x] \circ [z]$.
 \item $[x]*[y]=[x+y]=[y+x]=[y]*[x]$.
 \item $\psi [z] = [z] \otimes [z]$.
 \item The sub-Hopf algebra of $F_{*}\underline{E}_{\,*}$ generated
  by all $[x]$ with $x \in G^{*}$ is the ring-ring of $G^{*}$ over
  $F_{*}$, i.e. $F_{*}[G^{*}]$.
\end{enumerate}

\section{Properties of $\hko{0}=$ $\hgy{(\ints \times BO)}$ }
\label{sec-hkozero}

In this section, we record the known mod 2 homology for $\ko{0}$ and
introduce Hopf ring properties for the elements in homology.  We will
compute 
$\hko{0}=\hgy{(\ints \times BO)}\cong \hgy{(\ints)}\otimes\hgy{(BO)}$,
where $\hgy{(\ints)}$ is concentrated in $\deg 0$ and $\ints$ is the
set of integers with the discrete topology.

To understand $\hgy{(BO)}$, we examine
 $$\rls P^{\infty}=1 \times BO(1)
   \subset 1 \times BO \subset \ints \times BO =\ko{0}.$$ 
Recall that the homology of real projective space is given by 
 $$\hgy{(\rls P^\infty)}=\ints/2\{b_i:i\geq 0\},$$
the free $\ints/2$-module on the elements $b_{i}$ with $\deg(b_i)=i$.
The map 
        $$\rls P^\infty\rightarrow 1\times BO,$$  
which classifies the unreduced canonical line bundle, induces an embedding 
in homology. We denote the image of $b_{i}$ under this map as 
$z_{i}\in \hgy{(1\times BO)}$.  The classical result is that
        $$\hgy{(1\times BO)}=P(z_{i}:i>0),$$
and therefore
                        $$\hgy{(BO)}=P(\olz{i}:i>0).$$

We may gain information about the elements in $BO$ by considering their 
corresponding behavior in $\rls P^\infty$:
\begin{enumerate}
 \item The coproduct in $\hgy{(\rls P^{\infty})}$ is given by
 $\psi(b_i)=\sum_{j+k=i}b_j\otimes b_k$, so the same holds true in
 $\hgy{(BO)}$:
         $$\psi(z_i)=\sum_{j+k=i}z_j\otimes z_k.$$
 \item If we define $z(a)=\sum z_{i}a^{i}$, we see that consistent with the 
 coproduct is the equality
         $$z(a)\circ z(b)=z(a+b).$$ 
 By comparing coefficients, we obtain 
 $$z_j\circ{z_k}=\frac{(j+k)!}{j!k!}z_{j+k}.$$
\end{enumerate}

Looking to $\ints$, we note that its mod 2 homology is given by the
ring-ring $\ints/2[\ints]$, whose polynomial generators are of the
form $[b]$, for $ b\in \ints$. The Hopf ring properties for these
elements can be found in Section 3. Note that there are two
equivalent ways to express the same element, $1=[0_0]$, in
$H_{0}\ko{0}$. Additionally, $z_0$ is the basis element for
$H_0(\ints\times BO)= \ints/2 [\pi_0(\ints\times BO)]$, corresponding
to the component $1\times BO$, and will also be denoted by the
ring-ring element $[1]$.  Accordingly, $\ints/2[\ints]$ can be written
as $P([1],[-1])=P(z_0^{\pm 1})$.

We denote the element $z_i*[-1]\in\hko{0}$ by $\olz{i}=z_i/z_0$, allowing us
to prove:
\begin{lem}
For any element $[x]$, 
$$[x]\circ\olz{i}=0 \ \mbox{iff} \ [x]\circ z_i=0.$$ \label{lem:xcircz}
\end{lem}
\begin{proof}
Suppose that $[x]\circ\olz{i}=0$. Then the distributive property shows                  $$[x]\circ
z_i=[x]\circ(\olz{i}*[1])=([x]\circ\olz{i})*
        ([x]\circ[1])=0.$$ 
Conversely, suppose $[x]\circ z_i=0$. Then
$$[x]\circ\olz{i}=[x]\circ(z_i*[-1])=([x]\circ z_i)*
        ([x]\circ[-1])=0.$$ 
\end{proof}
                          
In summary, our Hopf ring relations thus far are given by
\begin{enumerate}
        \item $z_0=[1]$ 
        \item $1=[0_0]$
        \item $ {\psi(z_i)=\sum_{j+k=i}z_j\otimes z_k}$ and 
        $\psi(z(t))=z(t)\otimes z(t)$
        \item $ {z_j\circ{z_k}=\frac{(j+k)!}{j!k!}z_{j+k}}$
        \item $\psi[z]=[z]\otimes[z]$
        \item $[x]*[y]=[x+y]$
        \item $[z]\circ[x]=[zx]$
\end{enumerate}  
and we  will begin our calculations with
        $$\hko{0}=\hgy{(\ints\times BO)}=\ints/2[\ints]\otimes 
        P(\olz{i}:i>0)=P(z_{0}^{\pm 1})\otimes P(\olz{i}:i>0)=$$
        $$P(\olz{i},[-1]).$$      

\section{Machinery}
                   
Before delving  into our calculations, we quickly review the tools used
throughout the rest of this article, namely the \bss \ and the Frobenius 
and Verschiebung maps.
             
\subsection{The Bar Spectral Sequence}

Let $E$ be an $\Omega$-spectrum, and let $\spec{n}'$ be the component
of the basepoint in the $n$'th space of $E$.  The \bss\ (also called
the Rothenberg-Steenrod spectral sequence) is really the homology
Eilenberg-Moore spectral sequence for the path-loop fibration on the
space $\spec{n+1}'$,
\[
\begin{diagram}
     \node{\Omega\spec{n+1}'=\Omega\spec{n+1}}
     \arrow{e}\arrow{s,l}{\cong}
     \node{P\spec{n+1}'}
     \arrow{e}\arrow{s,l}{\cong}\node{\spec{n+1}'.}\\
     \node{\spec{n}} \node{\star}
\end{diagram}
\]

Based on the bar resolution, it is a spectral sequence of Hopf 
algebras with the basic property that
        $$ E^{2}_{s,t}= \mbox{Tor}^{H_{*}(\ul{E}_{\,n})}_{s,t}  
        ({\bf F},{\bf F}) \Rightarrow H_{s+t}(\spec{n+1}';{\bf F}),$$
where ${\bf F}=\ints/2$.  Its differentials are given by
        $$d_r:E^r_{s,t} \rightarrow E^r_{s-r,t+r-1}.$$  
For more information, see \cite{mccleary}.

When $R$ is a bicommutative Hopf algebra over ${\bf F}=\ints/2$, then
$\tor{R}=\mbox{Tor}^R_{**}({\bf F},{\bf F})$ is again a bicommutative
Hopf algebra over ${\bf F}$.  If $I$ is the augmentation ideal, there
is a natural isomorphism 
$\si:I/I^2\rightarrow\mbox{Tor}^R_{1,*}({\bf F},{\bf F})$.  Moreover,
if $y\in I_k$ and $y^2=0$ there are naturally defined divided powers
$\gamma_j(\si(y))\in\mbox{Tor}^R_{j,jk}$.  The following computations for
$\mbox{Tor}^{R}_{*,*}(\ints/2,\ints/2)$ are well known (see
\cite{ravwil2}):
\begin{enumerate}
 \item  $\tor{P(x)} = E[\si(x)]$, the exterior algebra on 
 $\si(x)\in \mbox{Tor}_{1,*}$.
 \item  $\tor{E(x)}=\Gamma[\si(x)]$, the divided power 
 algebra on $\si(x)$, which is spanned by the elements $\gamma_j(\si(y))$.
 \item  $ \tor{{\bf Z}/2[{\bf Z}]}=E[\si(x)]$, where  
 $\ints/2[\ints]=\ints/2(x, x^{-1})$ is the ring-ring.
 \item  $ \mbox{Tor}^{A}\otimes \mbox{Tor}^{B} \stackrel{\simeq}
 {\longrightarrow}\mbox{Tor}^{A \otimes B}.$
\end{enumerate}

\subsection{Frobenius and Verschiebung Maps}

Given $A$, a bicommutative Hopf algebra over $\ints/p$ ($p$ prime), the 
Frobenius map $F:A\rightarrow A$ is described by $F(x)=x^{p}$. The 
Verschiebung map $V$ is defined to be the dual of $F$ on $A^{*}$, the dual 
of $A$.  In other words, if $F^\prime:A^*\rightarrow A^*$ is defined by $x 
\mapsto x^p$, then $(F^{\prime})^{*}=V: A \rightarrow A$.

{}From \cite{ravwil2}, these maps have the following properties (for
$p=2$):

\begin{enumerate}
 \item With a shift of grading, $V$ and $F$ are Hopf algebra maps.
 \item $VF=FV:$
         $$V(x^{2})=VF(x)=FV(x)=[V(x)]^{2}.$$
 \item $V(x\circ y)=V(x) \circ V(y)$
 \item $F(x \circ V(y))=F(x) \circ y$
 \item For the coalgebra $\Ga (x)$,
         $$V(\gamma_{2q}(x))=\gamma_{q}(x)$$
         $$V(\gamma_{q}(x))=0  \ \mbox{if}\ q\neq 0 \ mod \ 2\,  \ $$
 As a consequence, for the elements $z_{i}\in \hko{0}$, we have both     
 $V(z_{2i})=z_i$ and $V(z_{2i+1})=0.$ 
\end{enumerate} 

\section{The Calculation of the Hopf Ring for $\hko{*}$}
      
We are ready to proceed to the calculation.  Each of the following
sections contains the computation of the mod 2 homology for one of the
spaces in the spectrum $KO$.  As we continue, we will introduce
various Hopf ring relations; we defer the proofs of these relations
until the ends of their corresponding sections.
       
\subsection{$\hko{1}=\hgy{(U/O)}$\label{sec:one}}

We start by inputting $\hko{0}=\hgy{(\ints\times BO)}=P(z_{0}^{\pm 
1})\otimes P(\olz{i}:i>0)$ into the \bss.
        $$\etwo =\torh{({\bf Z}\times BO)}=\tor{P(z_0^{\pm 1})\otimes   
        P(\olz{i}:i>0)}$$ 
        $$=E(\si(z_0))\otimes E(\si(\olz{i}:i>0))=E(\si(z_i))\RA \hko{1}=
    \hgy{(U/O)}.$$
Since the elements $\si(z_i)$ are all in the first filtration, the \bss 
\ collapses at the $E^2$-term, and the $E^\infty$-term is $E(\si(z_i))$.

For any $x \in H_{*}\underline{E}_{\,n}$, \mbox{$e \circ x \in H_{*} 
\ul{E}_{\,n+1}$} detects (the image in $E^{\infty}$ of) $\si(x)$, 
where $e$ is the fundamental class in $H_1\ko{1}$. As such, the element 
$e\circ z_i$ detects $\si(z_i)$. The coproduct for $e$ is given by 
        $$\psi e=e\otimes 1+ 1\otimes e=e\otimes[0]+ [0]\otimes e.$$

It is useful to note the following:
\begin{lem}
For all elements $a$ with $\varepsilon a=0$, 
                $$e\circ(a*[1])=e\circ a.$$
In particular, for any $i>0$, $e\circ\olz{i}=e\circ z_i$. \label{lem:ecircz}
\end{lem}
\begin{proof} Using distributivity we obtain 

\begin{tabular}{r l } 

$e\circ(a*[1]) =$ & $(e\circ a)*(1\circ [1])+(1\circ a)*(e\circ [1])$\\
        $=$ & $(e\circ a)*[0]+0*(e\circ[1])=e\circ a$. \\
\end{tabular}

Thus $e\circ z_i=e\circ(\olz{i}*[1])$ is equivalent to $e\circ\olz{i}$
for all $i>0$.
\end{proof}

Returning to our calculation, we wish to simplify the
$E^{\infty}$-term as a Hopf ring.  To do so, we use the relation
\begin{equation}
 F(\es{n})=(e^{\circ n})^2=e^{\circ{n}}\circ{z_n}\label{eq:etonsq}
\end{equation}
(to be proved below) with $n=1$. Since $V(z_{2i})=z_i$ and $F(x\circ
V(y))=F(x) \circ y$, we may write
        $$(e\circ z_i)^2=F(e\circ z_i)=F(e)\circ z_{2i}= (e \circ 
        z_1)\circ       z_{2i}=e \circ z_{2i+1}.$$
Note that we have applied the $\circ$-product relation for the elements 
$z_{i}$, from section \ref{sec-hkozero}.

We apply this calculation as often as possible, by using the fact that any 
positive integer $m$ can be expanded in its binary form; 
$m=2^q(2i+1)-1=1+2+2^{2}+\ldots+2^{q-1}+2^{q+1}i$ with $q\geq 0, i\geq 0$, 
both integers.  Hence, every element $e\circ z_m$ is equivalent to  the 
product $F^q(e \circ z_{2i})$, and so each element $e\circ z_{2i}$ is a 
polynomial generator.

Consequently,  $\hko{1}=\hgy{(U/O)}=P(e\circ z_{2i})$.

\subsubsection*{Proofs of Relations for $\hko{1}$}

{\bf {\it Relation (\ref{eq:etonsq}): $F(\es{n})=(\es{n})^2=\es{n}\circ 
z_n.$}} 

\begin{proof} Let $X$ be an infinite loop space.  For each $n\geq 0$, the 
Dyer-Lashof operations $Q^n:H_q(X)\rightarrow H_{q+n}(X)$ are defined.  
For our purposes, we need the following properties (from \cite{curtis} and 
\cite{priddy}):
\begin{enumerate}
        \item  $Q^n([1])={z_n}*[1]$
        \item  $Q^n(u)=u^2$ if dim$(u)=n$.
        \item  The homology suspension homomorphism $\si:H_*(\Omega X)
        \rightarrow H_{(*+1)}(X)$ yields $\si(Q^n(u))=Q^n(\si(u))$.
\end{enumerate}
In our case, $\si(x)=e\circ x$, so
        $$F(\es{n})=(\es{n})^{2}=Q^n(\es{n})=Q^n(\es{n}\circ[1])=
        \es{n}\circ Q^n([1])=\es{n}\circ(z_n*[1]).$$

By lemma \ref{lem:ecircz}, $\es{n}\circ(z_n*[1])=\es{n}\circ{z_n}. $
\end{proof}

\subsection{ $\hko{2}=\hgy{(Sp/U)}$}

The \bss \ produces 
        $$\etwo = \torh{(U/O)}=\tor{P(e\circ z_{2i})}=E(\si(e\circ 
        z_{2i}))$$ 
        $$\RA \hgy{(Sp/U)}.$$ 
The elements $\si(e\circ z_{2i})$ are all in the first filtration, 
collapsing the \bss \ at the $E^2$-term. Accordingly, the $E^\infty$-term 
is \mbox{$E(\si(e\circ z_{2i}))$}.


The rest of the proof is the same as the proof of 6.1, except we use property 3 of page 3 instead of property 2.  Then we have $\hko{2}=\hgy{(Sp/U)}=P(e^{\circ 2}\circ z_{4i})$.

\subsection{ $\hko{3}=\hgy{(Sp)}$}

The \bss \ gives 
        $$\etwo = \torh{(Sp/U)}= \tor{P(e^{\circ 2}\circ z_{4i})}=
    E(\si(e^{\circ 2}\circ z_{4i}))$$ 
    $$\RA \hgy{(Sp)}$$
Again, each generator is located in the first filtration.  Accordingly, 
the \bss \ collapses at the $E^2$-term, and the $E^\infty$-term is given 
by $E(\si(e^{\circ 2}\circ z_{4i}))$.

%


The rest of the proof is the same as the proof of 6.1, except we use the property
\begin{equation}
       (\es{3})^2=0\label{eq:etonsqu}
\end{equation}
instead of property 2.  We have $\hko{3}=\hgy{(Sp)}= E(e^{\circ 3}\circ z_{4i})$.

\subsubsection*{Proofs of Relations for $\hko{3}$}

{\it Relation (\ref{eq:etonsqu}): $(\es{3})^2=0.$}

\begin{proof}
For any elements $a, b$  with $\varepsilon a=\varepsilon b=0$, we may use 
the distributive property to obtain $e\circ(a*b) 
=(e\circ{a})*(1\circ{b})+(1\circ{a})*(e\circ{b})=0$, using the coproduct 
for $e$.  Thus $e\circ(\mbox{decomposables})=0.$ We use this fact and 
consequences of relation (\ref{eq:etonsq}) with $n=1, 3$ to compute 
        $$(e^{\circ{3}})^{2}=e^{\circ{3}}\circ{z_{3}}= \es{2}\circ e
        \circ{z_{1}}\circ{z_{2}}=\es{2}\circ(e)^{2}\circ{z_{2}}=0. $$ 
\end{proof}

\subsection{ $\hko{4}=\hgy{(\ints \times BSp)}$}

Since the connected part of $\ko{4}$ is $BSp$, the \bss \ gives 
        $$\etwo = \torh{(Sp)}= \tor{E(e^{\circ 3}\circ z_{4i})} = 
        \Ga(\si(e^{\circ 3}\circ z_{4i}))$$ 
        $$\RA \hgy{(BSp)}.$$ 
The \bss \ collapses at the $E^{2}$-term, as each element has even total 
degree. The $E^{\infty}$-term is therefore given by $\Ga(\si(\es{3}\circ 
z_{4i}))$, or equivalently by $E(\ga(\si(e^{\circ 3}\circ z_{4i})))$. 

We will simplify in stages, starting with the elements  in the first 
filtration, that is, the elements $\si(\es{3}\circ z_{4i})$. As usual, 
$e^{\circ 4}\circ z_{4i}$ detects $\si (e^{\circ 3} \circ z_{4i})$. 
We use relation (\ref{eq:etonsq}) with $n=4$ and the equality 
$V(z_{8i})=z_{4i}$ to determine the algebra structure:
        $$(e^{\circ 4}\circ z_{4i})^2=F(e^{\circ 4}\circ z_{4i})=
    F(e^{\circ 4})\circ z_{8i}=(e^{\circ 4} \circ z_4)\circ z_{8i}=
    e^{\circ 4} \circ z_{8i+4}.$$  
Again, expanding every positive integer $m$ in its binary form 
\mbox{$2^q(2i+1)-1$} yields $\es{4}\circ z_{4m}=F^q(\es{4}\circ z_{8i})$.  
Thus the elements in the first filtration are polynomial generators given by 
$\es{4}\circ z_{8i}$.

Next, we examine the elements of $\ga, j>0$.  As this calculation will make repeated use of the 
Verschiebung map, we record
the following:

\begin{prop}The Verschiebung map has the properties
        \begin{enumerate}
\item the element $V^j(\ga(\si(z_i)))$ is detected by the element $e\circ z_i$.         
\item if $V^j(x)=z_i$, then $x=z_{2^j \cdot i}$, mod decomposable elements.
        \end{enumerate}
\label{prop:vs}
\end{prop}

\begin{proof} 
\begin{enumerate}
\item By virtue of the fact that $V(\gamma_{2q}(x))=\gamma_q(x)$, we obtain
$V^j(\ga(\si(z_i)))=V^{j-1}(\gamma_{j-1}(\si(z_i)) =...=\si(z_i)=e\circ z_i$.
\item The facts that $V(z_{2i})=z_i$ and $V(z_{2i+1})=0$ may be used to prove this equation.     
\end{enumerate}
\end{proof}

Suppose $x$ detects $\ga(\si(\es{3}\circ 
z_{8i}))$. Then by proposition \ref{prop:vs},
        $$V^jx=\es{4}\circ z_{8i}.$$ 
Further, 
        $$V^jF^qx=F^qV^jx=F^{q}(\es{4}\circ z_{8i})=\es{4}\circ z_{4m}$$ 
detects $\si(\es{3}\circ z_{4m})$ and thus $F^q x$ detects 
$\ga(\si(\es{3}\circ z_{4m}))$.  Thus the elements $x$ (as $i$ and $j$ 
vary) are polynomial generators for $\hgy{(BSp)}$.  

All that remains is to identify the elements $x$. To do so, we map down 
from $\ko{4}$ to $\ko{-4}= \ko{4}$ via the map $\circ [\la]$, and seek instead 
to simplify $x\circ [\la]$.  Application of the relation 
\begin{equation}
        \es{4}\circ[\la]=\olz{4}\circ[\be]\label{eq:ecirclam} 
\end{equation} 
produces the following, modulo decomposable elements:  
        $$V^j(x\circ [\la])=\es{4}\circ z_{8i}\circ [\la]=\olz{4}\circ  
        z_{8i}\circ[\be]=\olz{8i+4}\circ[\be].$$ 
Accordingly, proposition \ref{prop:vs} shows that 
        $$x\circ[\la]=\olz{2^j(8i+4)}\circ[\be]+\mbox{decomposables.}$$ 
By virtue of the relations 
\begin{equation}
        \olz{1}\circ[\be]=z_{1}\circ[\be]=\olz{2i+1}\circ[\be]=0
        \label{eq:olzcircbe}
\end{equation} 
\begin{equation}
        \olz{2}\circ[\be]=z_{2}\circ[\be]=\olz{4i+2}\circ[\be]=0
        \label{eq:evenzcircbe}
\end{equation} 
and by the property that every positive number, divisible by 4, has 
the unique form $2^j(8i+4)$, we obtain $\hgy{(BSp)}=P(\olz{4i}\circ[\be]: 
i>0)\subset\hko{-4}.$ Thus
        $$\hgy{\ko{-4}}=\hgy{(\ints)} \otimes \hgy{(BSp)}=P([\be][-
    \be])\otimes P(\olz{4i}\circ[\be]:i>0)=$$ 
    $$P(\olz{4i}\circ[\be],[-\be]).$$ 

To find $\hko{4}$, we use $\circ [\la^{-1}]$ to map up 8 spaces. 
We obtain $\hko{4}=$ $\hgy{(\ints \times BSp)} =$ $ P(\olz{4i}
\circ[\be\lai],[-\be\lai]).$

\subsubsection*{Proofs of Relations for $\hko{4}$}

{\it Relation (\ref{eq:ecirclam}): $\es{4}\circ[\la] 
=\olz{4}\circ[\be].$}

\begin{proof}
There is only one nontrivial element in the bidegree of $\olz{4}\circ[\be]$ so either this is zero or it is as claimed above.  If it is zero then if we take circle multiplication by $\es{4}$ it should still be zero. However,

\begin{tabular} {r l  l}
        $\es{4}\circ[\be]\circ\olz{4}=$ & $\es{2}\circ\es{2}\circ[\be]\circ\olz{4}$ &\\
                                =& $\es{2}\circ[\eta^2]\circ\olz{2}\circ\olz{4}$ & by 7 page 3\\
                                =& $(e\circ[\eta])^{\circ 2}\circ\olz{2}\circ\olz{4}$ & \\
                =& $(\olz{1})^{\circ 2}\circ\olz{2}\circ\olz{4}$ & by 6 page 3\\
= & $\olz{1}\circ \olz{7}$& by 1 page 3\\
= & $\olz{1}^2\circ \olz{6}$ & by 9 page 3\\
= & $(\olz{1}\circ \olz{3})^2$& by distributivity\\
= & $((\olz{1})^2\circ\olz{2})^2$ & by 9 page 3\\
= & $(\olz{1}\circ\olz{1})^4$ & by distributivity\\
= & $(\olz{1})^8$ & by 9 page 3\\
\end{tabular}

\noindent Which we know to be non-zero. Thus the relation has to hold.

\end{proof}

\noindent {\it Relation (\ref{eq:olzcircbe}): 
$\olz{1}\circ[\be]=z_{1}\circ[\be]= \olz{2i+1}\circ[\be]=0.$}

\begin{proof}
We first prove 
\begin{equation}
 e\circ[\eta]=\olz{1}.\label{eq:ecircal}
\end{equation}
By definition of $\eta\in\pi_1(BO)$, the homotopy element
$\eta:S^1\rightarrow BO$ classifies the reduced canonical line bundle.
Thus $e\circ[\eta]\in H_1(BO)$ is the image under $\eta$ of the
fundamental class $e$, which is nonzero by the Hurewicz theorem.  The
only possible element in this dimension (in $0\times BO$) is
$\olz{1}$, proving relation (\ref{eq:ecircal}).

We return to the proof of relation (\ref{eq:olzcircbe}). Using
relation (\ref{eq:ecircal}) and the fact that $\eta\be=0$ gives
$$ \olz{1}\circ [\be]=e\circ[\eta]\circ[\be]=
   e\circ[\eta\be]=e\circ[0]=0. 
$$ 
By lemma \ref{lem:xcircz},
        $z_{1}\circ[\be]=0. $
The distributive property completes the proof:
        $$\olz{2i+1}\circ[\be]=(z_{2i+1}*[-1])\circ[\be]=(z_{2i+1}\circ         
        [\be])*[-\be]=(z_{2i}\circ z_{1}\circ[\be]) *[-\be]=0. $$ 
\end{proof}

\noindent {\it Relation (\ref{eq:evenzcircbe}): 
$\olz{2}\circ[\be]=z_{2} \circ[\be]=\olz{4i+2}\circ[\be]=0.$}

\begin{proof} We refer to $KU$, the unitary Bott spectrum. As we will prove in section 8,
  $$H_{*}\ul{KU}_{\,0}= H_{*}(\ints\times BU)=P(\olz{2i},[-1]).$$
 We start by proving the relation
 \begin{equation}
   e^{\circ 2}\circ[\nu]=\olz{2},\label{eq:compl}
\end{equation}
where $e$ is the fundamental class in $\hku{1}$.  Since
$KU_{*}=\ints[\nu^{\pm 1}],$ with $\deg(\nu)=2$, relation
(\ref{eq:compl}) is proved using the same argument as relation
(\ref{eq:ecircal}) above.

Relation (\ref{eq:evenzcircbe}) may now be proven using the
forgetful map, $f:KU\rightarrow  KO$.  We build on 
relation (\ref{eq:compl}) to yield 
$e^{\circ 2}\circ[\nu]\circ[\nu^{2}]=e^{\circ 2}\circ[\nu^{3}]=
 \olz{2}\circ[\nu^{2}]$.   
Applying the forgetful map gives
 $$ 0 = \es{2}\circ[0]=f(e^{\circ 2}\circ[\nu^{3}])    
      = f(\olz{2} \circ[\nu^{2}]) = \olz{2}\circ[\be].$$  
In much the same way as the proof of relation (\ref{eq:olzcircbe}), we may 
use Hopf ring relations to obtain  the remaining equalities
 $$z_{2}\circ[\be]=0 \ \mbox{and} \ \olz{4i+2}\circ[\be]=0. $$
\end{proof}

\subsection{ $\hko{5}=\hgy{(U/Sp)}$}

The \bss \ gives 
 $$\etwo = \torh{({\bf Z}\times BSp)}= \tor{P(\olz{4i}\circ[\be\lai],
    [-\be\lai])}=$$
 $$E(\si(\olz{4i}\circ[\be\lai])) \RA \hgy{(U/Sp)}.$$ 
Due to the fact that every generator is once again located in the first 
filtration, the \bss\ collapses at the $E^2$-term. The $E^\infty$-term is 
thus given by $E(\si(\olz{4i}\circ[\be\lai]))$.

The element $e\circ \olz{4i}\circ[\be\lai]$ detects 
$\si(\olz{4i}\circ[\be\lai])$. Utilization of relation 
(\ref{eq:olzcircbe}) and lemma \ref{lem:ecircz} 
yields 
        $$F(e\circ \olz{4i}\circ[\be\lai])=F(e\circ z_{4i}\circ[\be\lai])=
    F(e)\circ z_{8i}\circ[\be\lai]=$$
    $$(e\circ z_1)\circ z_{8i}\circ [\be]\circ [\lai]=0.$$

As there are no further simplifications, $\hko{5}=\hgy{(U/Sp)}=E(e\circ 
z_{4i}\circ[\be\lai]).$

\subsection{ $\hko{6}=\hgy{(O/U)}$}

Since the connected portion of $\ko{6}$ is $SO/U$, the \bss \ produces
        $$\etwo = \torh{(U/Sp)}= \tor{E(e\circ z_{4i}\circ[\be\lai])}=$$ 
    $$\Ga(\si(e\circ z_{4i}\circ[\be\lai]))\RA \hgy{(SO/U)}.$$  
We express $\Ga(\si(e\circ z_{4i}\circ[\be\lai]))$ in its equivalent form 
$E(\ga(\si(e\circ z_{4i}\circ[\be\lai])))$, and note that since every 
element is located in even total degree, the \bss \ collapses at the 
$E^2$-term. Thus the $E^\infty$-term is given by $E(\ga(\si(e\circ 
z_{4i}\circ[\be\lai])))$.

As in $\hko{4}$, we simplify in stages, starting with elements of the first 
filtration, $\si(e\circ z_{4i}\circ[\be\lai])$. As usual, $\es{2}\circ 
z_{4i} \circ[\be\lai]$ detects $\si(e\circ z_{4i}\circ[\be\lai])$. The relation 
\begin{equation}
        \es{2}\circ [\be]=\olz{2}\circ [\eta^2] \label{eq:escircbe}
\end{equation}  
allows us to equate 
        $$\es{2}\circ z_{4i}\circ[\be\lai]=\olz{2}\circ z_{4i}
        \circ[\eta^2\lai]=\olz{4i+2}\circ[\eta^2\lai],$$ 
modulo decomposable elements. Therefore, in the first filtration we have 
an exterior algebra with generators $\olz{4i+2}\circ[\eta^2\lai]$.

Next, we proceed to the simplification of elements in $\ga$. Let $x$ 
detect the exterior algebra generator $\ga(\si(e\circ z_{4i}\circ
[\be\lai])).$ Then by proposition \ref{prop:vs},
        $$V^jx=\es{2}\circ z_{4i}\circ[\be\lai]=\olz{4i+2}\circ[\eta^2\lai]+    
        \mbox{decomposables}.$$ 
Thus, modulo decomposables, 
    $$x=\olz{2^j(4i+2)}\circ [\eta^2\lai].$$ 
Since every even positive number can be uniquely written as $2^j(4i+2)$, 
we must have $x=\olz{2i+2}\circ [\eta^2\lai]+\mbox{decomposables}$. Since 
\begin{equation} 
    \olz{2i+1}\circ [\eta^2]=0, \label{eq:olzcircal}
\end{equation}
we have $\hgy{(SO/U)}=E(\olz{2i+2}\circ[\eta^2\lai])$.

Thus $\hko{6}=\hgy{(O/U)}=\hgy{(SO/U)}\otimes \ints[(\ints/2)\eta^2\lai]= 
E(\olz{2i}\circ[\eta^2\lai])$

\subsubsection*{Proofs of Relations for $\hko{6}$}

{\it Relation (\ref{eq:escircbe}): 
$\es{2}\circ[\be]=\olz{2}\circ [\eta^2]$.}

\begin{proof} We use the forgetful map $f_{U}:KU\rightarrow  KO$, 
again building on relation (\ref{eq:compl}) to yield
        $$\es{2}\circ[\be]=f_{U}(e^{\circ 2}_{\bf C}\circ[\nu]\circ[\nu])       
        =f_{U}(\overline{z}_{{\bf C} ,2}\circ[\nu])=\olz{2}\circ[\eta^{2}]. $$
\end{proof}
      
\noindent {\it Relation (\ref{eq:olzcircal}): $\olz{2i+1}\circ [\eta^2]=0.$}

\begin{proof}
By relation (\ref{eq:ecircal}) and by the fact that  $\eta^{3}=0$ we have
        $$\olz{1}\circ[\eta^{2}]=e\circ[\eta]\circ[\eta^{2}]=e\circ[\eta^{3}]
        =e\circ [0]=0.$$ 
The argument to show $\olz{2i+1}\circ[\eta^{2}]=0$ is much the same as in 
relation (\ref{eq:olzcircbe}).
\end{proof}

\subsection{ $\hko{7}=\hgy{(O)}$}

The connected portion of $\ko{7}$ is $SO$, so the \bss \ gives 
        $$\etwo = \torh{(O/U)}= \tor{E(\olz{2i}\circ[\eta^2\lai])}=$$
        $$\Ga(\si(\olz{2i}\circ[\eta^2\lai])\RA \hgy{(SO)}$$ 
We rewrite $\Ga(\si(\olz{2i}\circ[\eta^2\lai]))$ as 
$E(\ga(\si(\olz{2i} \circ[\eta^2\lai])))$.

The elements $\ga(\si(\olz{2i}\circ [\eta^2\lai]))\in E^2_{2^j, 2^j(2i)}$ 
have total degree $2^j(2i+1)$. The elements 
        $$\gamma_1(\si(\olz{2i}\circ[\eta^2\lai]))=\si(\olz{2i}
        \circ[\eta^2\lai])$$ 
must survive to $E^\infty$, since any differential which
originates at one of these elements must end below the $x$-axis.  The images of  
these elements in $E^\infty$ represent the elements $e\circ\olz{2i}
\circ[\eta^2\lai].$

We now compare our spectral sequence with the bar spectral sequence for 
${\ints}/2$. Define the reflection map $r^{\prime}:B{\ints}/2=\rls P^{\infty} \rightarrow O$ 
as the map which sends a line $L\in \rls P^{\infty}$ to its reflection in 
the hyperplane orthogonal to $L$.  If we fix a standard line $L_0$, we 
may define the map $\displaystyle{r:\rls P^{\infty}
\rightarrow SO}$ by $\displaystyle{r(L)=r^{\prime}(L) r^{\prime}(L_0)}$.
Using $r$, we obtain the map $\Omega(r):\ints/2\rightarrow O/U$.  

We use $r$ to map from the spectral sequence
\begin{equation}
        \mbox{Tor}^{{\bf Z}/2[{\bf Z}/2]}(\ints/2, \ints/2) \Rightarrow \hgy{\rls 
        P^\infty}\label{eq:ss}
\end{equation} into our spectral sequence for $\ko{7}$.
The $\etwo$ term of the spectral sequence in (\ref{eq:ss}) is exterior 
on the elements $\ga(\si([\eta^2]))$.  Since there is no room for
differentials, the elements $\ga(\si([\eta^2]))$ survive in (\ref{eq:ss}), 
and therefore the elements $\ga(\si([\eta^2\lai]))$ also survive in our 
spectral sequence for $\ko{7}$.
             
Using the familiar properties of the Verschiebung map, we will show that 
the elements $\ga(\si(\olz{2i}\circ[\eta^2\lai]))$ must therefore also survive.
Since 
        $$\ga(\si(\olz{2i}\circ[\eta^2\lai]))=\ga\left(V^{j}(\olz{2^{j}(2i)}\circ 
        \si([\eta^{2}\lai]))\right)=\olz{2^j(2i)}\circ\ga(\si([\eta^2\lai])),$$
we may use the fact that the differentials respect the $\circ$-product, 
so all of our elements must survive as claimed.

Thus the $E^{\infty}$-term is given by $E(\ga(\si(\olz{2i} 
\circ[\eta^2\lai])))$.

As in $\hko{6}$, we simplify in steps.  The element $e\circ\olz{2i}\circ 
[\eta^2\lai] =e\circ z_{2i}\circ[\eta^2\lai]$ detects 
$\si([\eta^2\lai]\circ\olz{2i})$. Relation (\ref{eq:ecircal}) allows us to  
rewrite 
        $$e\circ z_{2i}\circ[\eta^2\lai]=\olz{1}\circ z_{2i}\circ[\eta\lai]=
        \olz{2i+1}\circ[\eta\lai],$$ 
modulo decomposable elements. 

Next, we let $x$ detect any exterior algebra generator 
$\ga(\si(z_{2i}\circ [\eta^{2}\lai]))$. Utilization of proposition 
\ref{prop:vs} shows 
        $$V^jx=e\circ z_{2i}\circ[\eta^2\lai]=\olz{2i+1}\circ[\eta\lai]+
    \mbox{decomposables}.$$ 
Hence, modulo decomposable elements, $x=\olz{2^j(2i+1)}\circ [\eta\lai].$ 
Since every positive number can be uniquely written as $2^j(2i+1)$, we 
must have $x=\olz{i+1}\circ [\eta\lai]+\mbox{decomposables}$. No further 
simplification is necessary, so $\hgy{(SO)}=E(\olz{i+1}\circ[\eta\lai]).$

We have therefore found $\hko{7}=\hgy{(O)}=\hgy{(SO)}\otimes 
\ints[(\ints/2)\eta\lai]= E(\olz{i}\circ[\eta\lai])$.

\subsection{ $\hko{8}=\hgy{(\ints\times BO)}$}

We complete the cycle, reaching $\ko{8}=\ints \times BO=\ko{0}$. The \bss 
\ gives 
        $$\etwo = \torh{(O)}= \tor{E(\olz{i}\circ[\eta\lai])}=$$ 
        $$\Ga(\si(\olz{i}\circ[\eta\lai]))\RA \hgy{(BO)}.$$  
We rewrite $\Ga(\si(\olz{i}\circ[\eta\lai]))$ in its equivalent form 
$E(\ga(\si(\olz{i}\circ[\eta\lai])))$.  The \bss \ collapses at the 
$E^2$-term, by comparison with the collapsing bar spectral sequence 
        $$\torh{(O(1))}=\ints/2\oplus \ints/2\{\gamma_{k+1}(\overline{[\eta]})\}$$
        $$\RA \hgy{(BO(1))}=\hgy{\rls P^{\infty}}=\ints/2\oplus         
        \ints/2\{\olz{k+1}\}.$$ 
The differentials for the original spectral sequence are therefore given 
by 
        $$d_r(\gamma_{k}\overline{[\eta]}\circ z(t))=d_r(\gamma_{k}\overline{[\eta]})\circ      
        z(t)=0,$$
where $z(t)=\sum_{i\geq 0}z_i t^i$.

The element $e\circ\olz{i}\circ[\eta\lai]$ detects $\si(\olz{i}\circ
[\eta\lai])$. Use of  relation (\ref{eq:ecircal}) allows us to rewrite 
$e\circ \olz{i}\circ[\eta\lai]$ as $\olz{1}\circ \olz{i}\circ[\lai].$ 

As usual, we start with elements in the first filtration.  Using the 
relation
\begin{equation}
        \olz{1}\circ\olz{2i+1}=\olz{1}^{2}\circ\olz{2i}.\label{eq:olzcircolz}
\end{equation}
gives
        $$\olz{1}\circ\olz{2i+1}\circ[\lai]=\olz{1}^{2}\circ\olz{2i}
        \circ[\lai]=F(\olz{1})\circ\olz{2i}\circ[\lai]=F(\olz{1}\circ\olz{i}
        \circ[\lai]).$$
We apply this calculation as often as possible, by expanding each integer 
in its binary form, $m=2^q(2i+1)-1$.  Then 
    $$\olz{1}\circ \olz{m}\circ[\lai]=F^q(\olz{1}\circ\olz{2i}\circ[\lai])
    =F^q(\olz{2i+1}\circ[\lai])+\mbox{decomposables}.$$ 
Therefore in the first filtration $\gamma_{1}$, we obtain a polynomial 
algebra with generators $\olz{2i+1}\circ[\lai]$.

We now look to simplify the remaining elements in $\ga$. By the preceding 
argument, we have a polynomial algebra with generators given by 
$x=\ga(\olz{2i+1}\circ[\lai])$. Proposition \ref{prop:vs} proves that since
$V^j(x)=\olz{2i+1}\circ[\lai]$, $x$ may be expressed as $\olz{2^j(2i+1)}\circ[\lai],$ mod decomposables.
Since every positive integer can be uniquely expressed in the form 
$2^j(2i+1)$, we have $\hgy{(BO)}=P(\olz{i+1}\circ[\lai])$.  

Thus $\hko{8}=\hgy{(\ints \times BO)}=P([\lai][-\lai])\otimes P(\olz{i+1} 
\circ[\lai]) =$ $P(\olz{i}\circ[\lai],[-\lai])$.
 
Since our answer for $\hko{8}$ is equivalent to our answer for $\hko{0}$, 
the Bott-periodic circle from $\ko{0}$ to $\ko{8}$ is complete.

\subsubsection*{Proofs of Relations for $\hko{8}$}

{\it Relation (\ref{eq:olzcircolz}):} $\olz{1}\circ 
\olz{2i+1}=\olz{1}^{2}\circ \olz{2i}.$

\begin{proof}
We first find the coproduct for $\olz{1}$;

\begin{tabular} {r l}
        $\psi(\olz{1})=$ & $\psi(z_{1}*[-1])=\psi(z_1)*\psi([-1])$\\
    $=$ & $((z_1\otimes z_0)+(z_0\otimes z_1))*([-1]\otimes [-1])$\\ 
    $=$ & $\olz{1}\otimes 1+1\otimes \olz{1}.$\\
\end{tabular}  

\noindent Thus $\olz{1}\circ \mbox{decomposables}=0$, as
$\olz{1}\circ(a*b)=(\olz{1}\circ a)*(1\circ b)+(1\circ 
a)*(\olz{1}\circ b)=0$.

For any $i>0$ we may use distributivity to simplify 
        $$\olz{1}\circ\olz{i}=\olz{1}\circ(z_i*[-1])=(\olz{1}\circ z_i)*(1 
        \circ [-1])+(1\circ z_i)*(\olz{1}\circ[-1])=\olz{1}\circ{z_i}+0.$$
Thus, for all $i>0$, $\olz{1}\circ\olz{i}=\olz{1}\circ{z_i}.$ In 
particular, $\olz{1}\circ\olz{2i}=\olz{1}\circ z_{2i}$.  Next,

\begin{tabular}{r l }
        $\olz{1}\circ\olz{2i}=$ & $\olz{1}\circ{z_{2i}}=(z_1*[-1])
        \circ{z_{2i}}$\\
        $=$ & $(z_1\circ{z_{2i}})*([-1]\circ{z_0})+ \mbox{      
                decomposables}$\\
        $=$ & $z_{2i+1}*[-1]+\mbox{decomposables}$\\
        $=$ & $\olz{2i+1}+\mbox{decomposables.}$\\
\end{tabular}

\noindent Since $\olz{1}\circ\olz{2i}=\olz{2i+1}+\mbox{decomposables}$, 
and since $\olz{1}\circ \mbox{decomposables}=0$, we now have
        $$\olz{1} \circ\olz{2i+1}=\olz{1}\circ\olz{1}\circ\olz{2i}      
        +\olz{1}\circ\mbox{ decomposables}=\olz{1}^{\, \circ 2}\circ    
        \olz{2i}.$$

To finish, we need only prove $\olz{1}^{\, \circ 2}=\olz{1}^{2}$.  We have 
already established that $\olz{1}\circ\olz{i}=\olz{1}\circ z_i$ for all 
$i>0$, and so  $\olz{1}\circ\olz{1}= \olz{1}\circ {z_{1}}$. Using 
distributivity and the $\circ$- product for the $z_{i}$ shows 

\begin{tabular}{r l}
        $\olz{1}\circ\olz{1}=$ & $\olz{1}\circ z_1=(z_1*[-1])\circ{z_1}$\\
        $=$ & $(z_1\circ{z_1})*([-1]\circ{z_0})+ (z_1\circ{z_0})*([-    
        1]\circ{z_1})$\\
        $=$ & $0+z_1*([-1]\circ{z_1})$\\
        $=$ & $z_{1}*(\chi z_{1}).$\\
\end{tabular}

\noindent To find $\chi z_1$, we calculate $\chi z_{0} =\chi[1]=[-1]
\circ[1]=[-1]$. Next, we use the Hopf ring property that 
$\sum a^{\prime}*\chi a^{\prime\prime}=\varepsilon a$. Since $\varepsilon 
z_{i}=\delta_{i0}$, we have 
        $$z_{1}*\chi z_{0}+z_{0}*\chi z_{1}=\varepsilon z_{1}=0.$$  
Thus    
        $$z_{1}*[-1]+[1]*\chi z_{1}=0.$$
We may now solve for 
        $$\chi z_{1}=z_{1}*[-2].$$ 
By virtue of these facts, we have obtained $\olz{1}\circ\olz{1}=z_{1} 
*(z_{1}*[-2])=\olz{1}^{2}$, completing our proof.  
\end{proof}

\section{Properties of $\hku{0}=$ $\hgy{(\ints \times BU)}$ }

In this section, we record the known mod 2 homology for $\ku{0}$ and 
introduce Hopf ring properties for the elements in homology.  We will 
compute $\hku{0}=\hgy{(\ints \times BU)}\cong \hgy{(\ints)}\otimes 
\hgy{(BU)}$, where $\hgy{(\ints)}$ is concentrated in $\deg 0$ and $\ints$ 
is the set of integers with the discrete topology.  

To understand $\hgy{(BU)}$, we examine the map
$$\rls P^\infty \rightarrow \cmpx P^{\infty}\rightarrow 1 \times BU
\subset \ints \times BU =\ku{0}.$$ Recall that the homology of $\rls
P^\infty$ is spanned by generators $b_i$.  We write $b_{2i}$ again for
the image of $b_{2i}$ in $H_*\cmpx P^\infty$ (the image of $b_{2i+1}$
is zero).  The elements $b_{2i}$ form a basis for $H_*\cmpx P^\infty$.
We write $z_{2i}$ for the image of $b_{2i}$ in $H_{2i}\ku{0}$, which
is also the image of our previous element $z_{2i}$ under the map
$m:\ko{0}\rightarrow\ku{0}$.  We also write
$\olz{2i}=z_{2i}/z_0=z_{2i}*[-1]\in H_{2i}(0\times BU)$.  The
classical result is that
        $$\hgy{(0\times BU)}=P(\olz{2i}:i>0),$$
and it follows that $\hgy{(\ints\times BU)}=P(z_{2i},[-1])$. 

As the map $m_*:H_*\ko{*}\rightarrow H_*\ku{*}$ preserves Hopf ring
structures, the elements $z_{2i}$ in $H_*\ku{*}$ have the same Hopf
ring properties as they did in $H_*\ko{*}$.

\section{The Calculation of the Hopf Ring for $\hku{*}$}

\subsection{$\hku{1}=\hgy{(U)}$\label{sec:kuone}}

We input 
$\hku{0}=\hgy{(\ints\times BU)}=P([1],[-1])\otimes P(z_{2i}:i>0)$
into the \bss.
 $$\etwo =
   \torh{({\bf Z}\times BU)} =
   E(\si(\olz{0}))\otimes E(\si(\olz{2i}:i>0))=
   E(\si(\olz{2i}))
 $$
 $$ \RA \hku{1}= \hgy{(U)}.
 $$ 
Since the elements $\si(\olz{2i})$ are all in the first filtration,
the \bss \ collapses at the $E^2$-term, and the $E^\infty$-term is
$E(\si(\olz{2i}))$.

The element $e\circ\olz{2i}$ detects $\si(\olz{2i})$.  From $H_*\ko{*}$
we have $e^2=e\circ z_1$ and $z_1$ maps to $0$ in $H_*\ku{*}$ so
$F(e)=e^2=0$ in $H_*\ku{*}$.  Thus 
$$ F(e\circ\olz{2i})=F(e\circ V(\olz{4i}))=F(e)\circ\olz{4i}=0.$$
This solves the extension problem and shows that 
$H_*U=E(e\circ z_{2i})$ (which is of course well-known by other
methods).

\subsection{$\hku{2}=\hgy{(\ints \times BU)}$\label{sec:kutwo}}

To finish the cycle of 2 spaces in $\hku{*}$, we examine $\hku{2}$.

The \bss \ gives
        $$\etwo =\torh{(U)}=\Ga(\si(e\circ\olz{2i}))\RA 
        \hgy{(BU)},$$
since the connected portion of $\ku{2}$ is $BU$.
The \bss \ collapses at the $E^{2}$-term, as each element has even total 
degree. The $E^{\infty}$-term is therefore given by 
$\Ga(\si(e\circ\olz{2i}))$, or equivalently by 
$E(\ga(\si(e\circ\olz{2i})))$. 

We map down from $\ku{2}$ to $\ku{0}$ via the map $\circ[\nu]$.  
The elements $\ga(\si(e\circ\olz{2i}))$
are detected by $\ga(\es{2}\circ\olz{2i})$.  Applying the map
$\circ[\nu]$ to these elements and using relation (\ref{eq:compl}) yields
$$\ga(\es{2}\circ\olz{2i}\circ[\nu])=\ga(\olz{2i}\circ\olz{2}).$$
A similar argument as in the proof of $\hko{8}$ shows that these 
elements form the polynomial algebra $P(\olz{2i+2})$.  Mapping back 
to $\ku{2}$ via the map $[\nu^{-1}]$ yields
$P(\olz{2i+2}\circ[\nu^{-1}])=\hgy{(BU)}\subset \ku{2}$. Thus
$$\hku{2}=\hgy{(\ints\times BU)}=P([\nu^{-1}],[-\nu^{-1}])\otimes
P(\olz{2i+2}\circ[\nu^{-1}])$$
$$=P(\olz{2i}\circ[\nui],[-\nui]),$$ 
completing our proof.

\end{document}